\documentclass[onecolumn]{article}

\usepackage[sort, numbers]{natbib}
\usepackage{url}
\usepackage{amsmath}
\usepackage{mathtools}
\usepackage{amssymb}
\usepackage{amsthm}
\theoremstyle{plain}
\usepackage{color}
\usepackage{here}

\newcommand{\red}[1]{\textcolor{black}{#1}}

\def\llongrightarrow{\relbar\joinrel\relbar\joinrel\relbar\joinrel\rightarrow}

\providecommand{\rarrow}[1]{\stackrel{#1}{\llongrightarrow}}

\everymath{\displaystyle}

\usepackage{xr}
\newcommand{\ba}{\begin{align}}
\newcommand{\ea}{\end{align}}

\def\E{\mathbb{E}}
\def\N{\mathbb{N}}

\providecommand\X[1]{\boldsymbol{X_{#1}}}
\providecommand\Z[1]{\boldsymbol{Z_{#1}}}

\providecommand\phib{\boldsymbol{\emptyset}}

\providecommand\phib{\boldsymbol{\emptyset}}
\providecommand\X[1]{\boldsymbol{S_{#1}}}

\newtheorem{theorem}{Theorem}

\newtheorem{proposition}[theorem]{Proposition}

\theoremstyle{plain}
\newtheorem*{objectives*}{Objectives}

\def\agedit#1{{\color{black}#1}}
\title{Antithetic Integral Feedback ensures robust perfect adaptation in noisy biomolecular networks}
\author{Corentin Briat$^\dag$, Ankit Gupta$^\dag$, Mustafa Khammash\thanks{Corresponding author (e-mail: mustafa.khammash@bsse.ethz.ch); $^\dag$Authors contributed equally.}\\
Department of Biosystems Science and Engineering (D-BSSE), \\
ETH--Z\"{u}rich, Basel, Switzerland}

\begin{document}

\maketitle


\begin{abstract}
Homeostasis is a running theme in biology. Often achieved through feedback regulation strategies, homeostasis allows living cells to control their internal environment as a means for surviving changing and unfavourable environments. While many endogenous homeostatic motifs have been studied in living cells, some other motifs may remain under-explored or even undiscovered. At the same time, known regulatory motifs have been mostly analyzed at the deterministic level, and the effect of noise on their regulatory function has received low attention. Here we lay the foundation for a regulation theory at the molecular level that explicitly takes into account the noisy nature of biochemical reactions and provides novel tools for the analysis and design of robust homeostatic circuits. Using these ideas, we propose a new regulation motif, which we refer to as {\em antithetic integral feedback}, and demonstrate its effectiveness as a strategy for generically regulating a wide class of reaction networks. By combining tools from probability and control theory, we show that  the proposed motif preserves the stability of the overall network, steers the population of any regulated species to a desired set point, and achieves robust perfect adaptation -- all with low prior knowledge of reaction rates. Moreover, our proposed regulatory motif can be implemented using a very small number of molecules and hence has a negligible metabolic load. Strikingly, the regulatory motif exploits stochastic noise, leading to enhanced regulation in scenarios where noise-free implementations result in dysregulation. Finally, we discuss the possible manifestation of the proposed antithetic integral feedback motif in endogenous biological circuits and its realization in synthetic circuits.
\end{abstract}

Perfect adaptation is that property of a biological system (e.g. a cell) which enables it to adapt to an external stimulus so that it maintains responsiveness to further stimuli. To be effective, such an adaptation mechanism must be robust, i.e. it must remain functional over a wide range of stimulus levels and system parameters. It was shown in \cite{Yi:00} that robust perfect adaptation in bacterial chemotaxis is achieved due to integral feedback control in the prevalent chemotaxis model \cite{Alon:99}. Other homeostatic systems have also been shown to realize integral feedback control. For example it was demonstrated in \cite{ElSamad:02} that calcium homeostasis in mammals relies on an integral feedback strategy to achieve perfect adaptation to persistent changes in plasma calcium clearance or influx, a property that enables mammals to maintain physiological levels of plasma calcium within tight tolerances in spite of varying demands for calcium.  In \cite{Muzzey:09}, integral feedback was implicated in the robust regulation of membrane turgor pressure in Saccharomyces cerevisiae. Following an osmotic shock, nuclear enrichment of the MAP kinase Hog1 adapts perfectly to changes in external osmolarity, a result of an integral feedback action that requires Hog1 kinase activity. \red{However, as some theoretical studies have suggested \cite{Ma:09,Shinar:10}, adaptation may not be solely related to integral control} .

\begin{figure}
    \centering
    \includegraphics[width=0.95\textwidth]{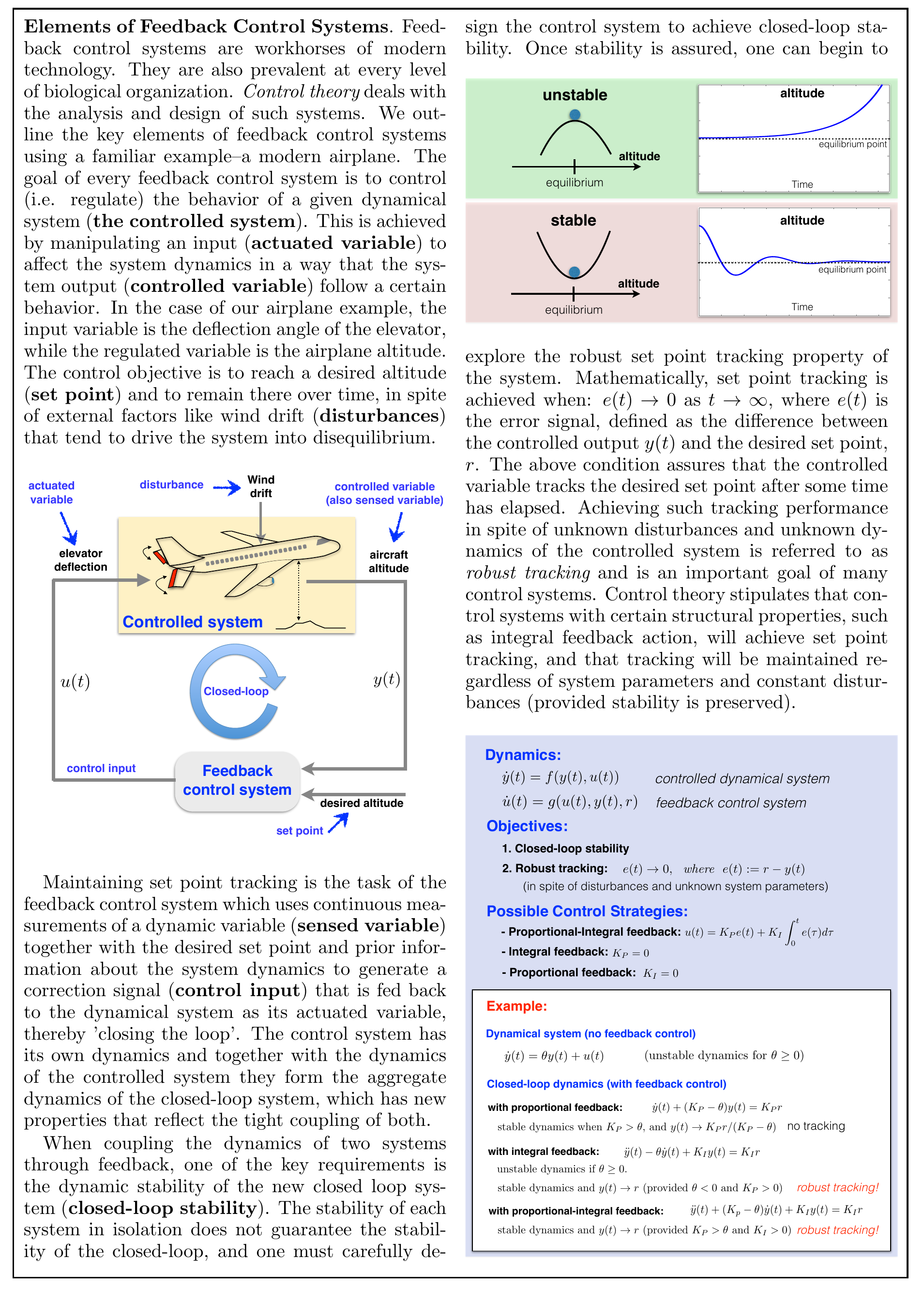}
\end{figure}

In engineering applications, integral feedback is recognized as a principal strategy for regulation. The Proportional-Integral-Derivative (PID) control architecture, which includes integral feedback as an essential element, is the workhorse of industrial control and is implemented in the majority of all automatic control applications \cite{Astrom:95}. Undoubtedly, the prevalence of such a control strategy in natural and man-made systems is due to the inherent property of integral feedback control to robustly steer a regulated system variable to a desired set point, while achieving perfect adaptation to disturbances (or stimuli), regardless of the  model parameters. Perhaps surprisingly, engineered biological circuits displaying perfect adaptation have received little attention so far, and current synthetic circuits only rely on simpler feedback strategies. For example, several control loops for controlling the level  of biofuel production in bacteria while still maintaining a low toxicity level are theoretically analyzed in \cite{Dunlop:10}. Another synthetic negative feedback loop is also designed in \cite{Stapleton:12} for the control of protein translation. Instead of integral feedback strategies, these circuits rely on the simpler  proportional feedback strategy. Consequently, they require a cumbersome tuning of parameters for achieving their goals. Such a tuning is very  difficult to realize in a biological setting, and even if a proportional feedback strategy is perfectly implemented, the absence of integral action implies that it does not possess the key property of perfect adaptation.

In the noise-free (deterministic) setting, integral feedback control is well-understood, and its ability to achieve robust set-point tracking and perfect adaptation is well-known \cite{Astrom:95}. In contrast, analogous strategies in intrinsically noisy cellular environments are unknown. Indeed in biologically important settings where the dynamics is described by stochastic processes (e.g. continuous-time discrete-state Markov processes), determining what constitutes integral feedback remains unclear. As in the deterministic case, a ``stochastic integral feedback" strategy must achieve closed-loop stability of the overall system, robust set-point tracking and robust perfect adaptation. Unlike the deterministic setting, however, set-point tracking and adaptation robustness must be maintained not only with respect to model parameters, but also for the highly fluctuating species abundances. One \red{possible way to construct} ``stochastic integral feedback" is to use statistical moments \red{(such as the mean or the variance)} to describe the process to be regulated, and then to design feedback regulation strategies that steer these moments to desired values while achieving perfect adaptation \cite{Khammash:11}. While this approach brings the problem back to the deterministic domain (statistical moments evolve according to deterministic dynamics), one is immediately faced with the moment closure problem, whereby an infinite set of differential equations is needed to determine even the first two moments; see e.g. \cite{Hespanha:08b}. Similar difficulties arise if one works with the chemical master equation; see e.g. \cite{Gillespie:92}.

\begin{figure}[H]
    \centering
    \includegraphics[width=0.5\textwidth]{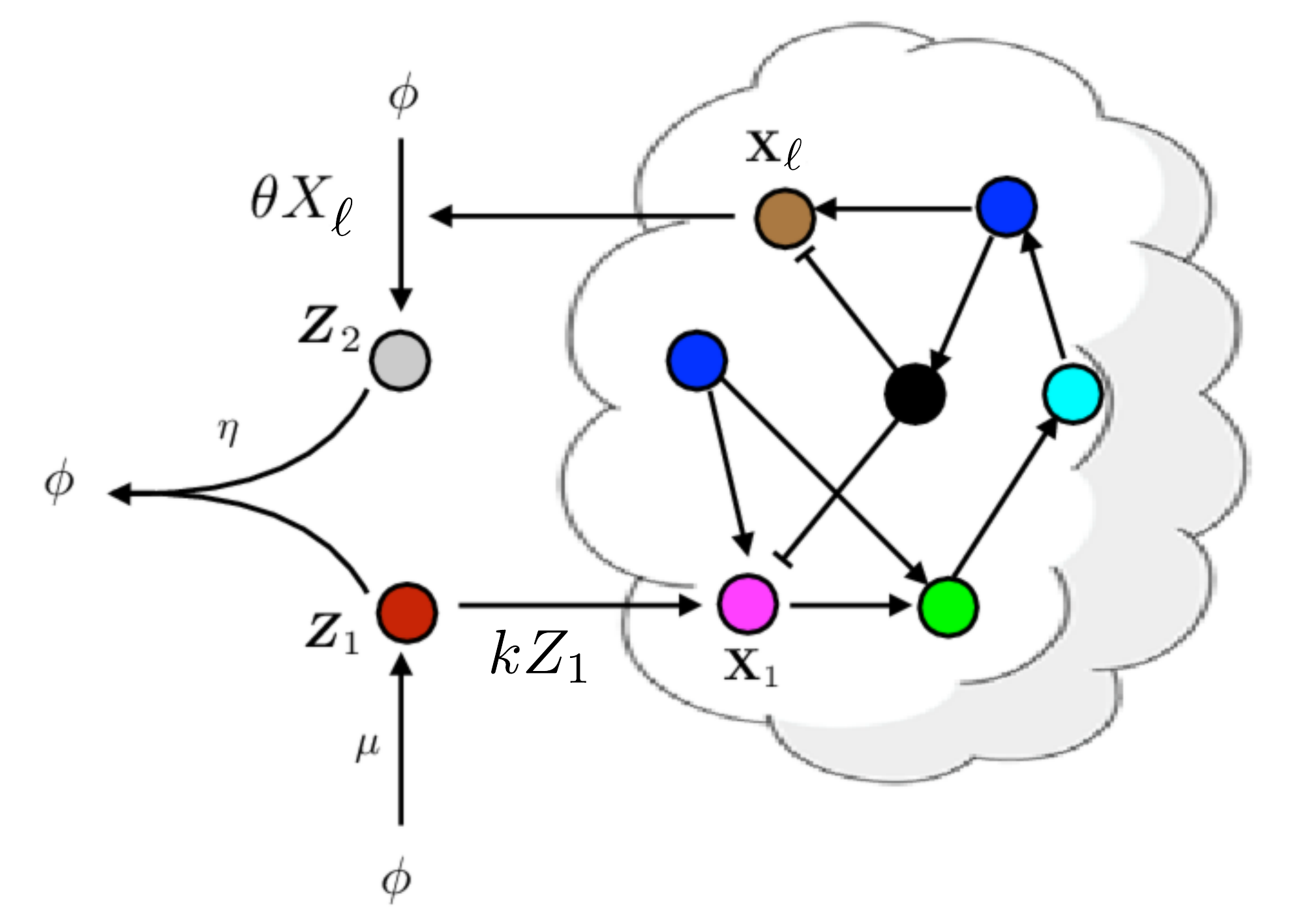}
 \caption{\red{Schematic representation showing the constituents of our biomolecular control system.  The network on the right (inside the cloud) represents the open-loop network, whose dynamics are to be controlled. Control is achieved by augmenting another network of reactions, referred to as the controller network (outside the cloud). Together the two networks form the closed-loop network whose dynamics are determined by the coupling resulting from the interaction of both networks.  The proposed controller network, which we refer to as antithetic integral controller,  acts on the open-loop network by influencing the rate of production of the actuated species $\X{1}$ by means of the control input species $\Z{1}$. The regulated species $\X{\ell}$ will be influenced by the increase or decrease of the actuated species $\X{1}$ and, in return, will influence the rate of production of the sensing species $\Z{2}$, that will, finally, annihilate with the control input species $\Z{1}$, thereby implementing a negative feedback control loop. The integral action is encoded in all the reactions of the controller network.}}\label{fig:cloud}
\end{figure}

Here we adopt a novel approach for designing a stochastic integral feedback strategy that exhibits robust \red{set-point} tracking and robust perfect adaptation. Rather than dealing with the deterministic moments dynamics, we work with the stochastic chemical reaction network directly, thereby circumventing the moment closure problem. The objective of our control setup, represented in Fig.~\ref{fig:cloud}, is to bring the population average of the species $\X{\ell}$ involved in a network (the ``cloud'' in Fig.~\ref{fig:cloud}) to a desired set-point. To achieve this, a new set of chemical reactions is introduced in a way that effectively implements a ``stochastic antithetic integral feedback controller". This controller network consists of four reactions and two additional controller species ($\Z{1}$ and $\Z{2}$) that can annihilate each other. The species $\Z{1}$ actuates the network which, in turn, influences the production of $\Z{2}$ through the output species $\X{\ell}$. We show in the results section that, for a large class of networks, the steady-state value for the population average of $\X{\ell}$ depends exclusively on \red{the ratio of} two of the controller parameters, and is independent of the network parameters. In this respect, the closed-loop network exhibits stability, robust set-point tracking, and robust perfect adaptation for $\X{\ell}$. To analyze such stochastic systems and to guarantee that they achieve these objectives, a new theory is needed. We develop such a theory here and use it to show that for a large class of networks, the considered antithetic feedback control motif can be used to achieve the desirable properties of ``stochastic integral feedback". We rigorously prove that such a motif robustly achieves the desired closed-loop stability (ergodicity) property. We additionally show that it achieves robust set-point tracking and robust perfect adaptation under mild conditions on the uncontrolled network. Intriguingly, our `` stochastic antithetic integral control motif" can provably achieve all the desired properties mentioned above, even when very low molecular copy numbers exist anywhere in the network. This presents a clear advantage in synthetic biology applications, where synthetic control loops involving large molecular counts can impose a debilitating metabolic load on the cell. Our control scheme can also be shown to possess remarkable stabilizing properties that are not found in deterministic implementations of the same circuit. This provides a clear example where the intrinsic stochastic noise is beneficial--it stabilizes a system which would otherwise be unstable. To the best of our knowledge, such a beneficial effect of noise, in the context of control theory, is reported here for the first time.  Note that many other benefits of noise, such as stochastic focusing \cite{Paulsson:00},  noise-induced oscillators \cite{Vilar:02} and noise-induced switches \cite{Tian:06,Acar:08}, have appeared in the literature in recent years.

\section*{Results}

In what follows, we elaborate on the control problem under consideration, the proposed controller, along with some technical results stating the conditions under which the proposed controller solves the considered control problem. Interestingly, these conditions obtained from probability theory connect to well-known concepts of control theory, such as stability and controllability. Some additional properties, such as robustness and innocuousness, are also discussed.

\subsection*{\red{The Network to be Controlled (Open-Loop Network)}}

\red{We start by describing} the reaction network we aim to control. \red{Consider a reaction network with mass-action kinetics involving $d$ molecular species denoted by $\X{1},\dots,\X{d}$. Under the well-mixedness assumption \cite{Anderson:11}, we can model the dynamics by a Markov process whose state at any given time is simply the vector of molecular counts of the $d$ species. The ``Markovian" assumption on the dynamics means that given the current state of the system, the future evolution of the state is independent of the past (memoryless property).} The state evolution is \red{influenced} by $K$ reaction channels: if the state is $x$, then the $k$-th reaction fires at rate $\lambda_k(x)$ and it displaces the state by the \emph{stoichiometric vector}, $\textstyle \zeta_k \in \mathbb{Z}^d$, \red{where $ \mathbb{Z}$ denotes the set of integers}. Here $\lambda_k$ is called the \emph{propensity function} of the $k$-th reaction and is assumed \red{to satisfy the property that} if for any $\textstyle x\in\mathbb{N}_0^d$  we have $\textstyle x+\zeta_k\notin\mathbb{N}_0^d$, then $\lambda_k(x)=0$, \red{where $\mathbb{N}_0$ denotes the nonnegative integers}. This property ensures that molecular counts of all the species remain nonnegative throughout the dynamics. \red{In the following, we shall refer to this network as the \emph{open-loop reaction network} and denote it by $(\X{},\lambda,\zeta)$.}

We now fix a state-space $\mathcal{S}$ for the Markovian reaction dynamics. This set $\mathcal{S}$ is a non-empty subset of $\textstyle \N^d_0$ which is closed under the reaction dynamics. This means that for any state $x \in \mathcal{S}$ we must also have $(x+\zeta_k) \in \mathcal{S}$ if the the $k$-th reaction has a positive rate of firing ($\lambda_k(x) >0$) \agedit{at state $x$}. Selecting $\mathcal{S}$ this way allows us to use it as a generic state-space for all Markov processes describing the reaction kinetics and starting at an initial state in $\mathcal{S}$ \cite{Briat:13i}. \red{Henceforth, we denote by} $\{X(t) = (X_1(t),\dots,X_d(t)) : t \geq 0\}$ the continuous time Markov process representing the reaction dynamics with an initial state $x_0 \in \mathcal{S}$.

From a control theoretic point of view, it is necessary to define input and output nodes of the above network. We assume here that species  $\X{1}$ is the \emph{actuated species} which is the species \red{that the controller} can act on. The \emph{regulated species} is $\X{\ell}$, for some $\ell\in\{1,\ldots,d\}$, and it is the species \red{we wish to control}. \red{The way \red{the controller acts on the actuated species, in order to control the regulated species} is depicted in Fig. \ref{fig:cloud}. This will be explained in more detail in the next section.}

\subsection*{\red{The Control Objectives}}  \red{We now state the objectives required from our control system, \agedit{which includes an open-loop network and a controller network.}}
\begin{objectives*}\label{problem}
 Find a controller \red{(set of additional reactions and additional species)} such that, by suitably acting on the actuated species $\X{1}$, we have the following properties for the closed-loop network (defined here as the interconnection of the open-loop network $(\X{},\lambda,\zeta)$ described above with the controller \agedit{network}):
  \begin{enumerate}
    \item the closed-loop network is ergodic;
      \item the first and second-order moments of $X(t)$ exist and are uniformly bounded and globally converging \agedit{with time} to their unique stationary value;
    \item we have that $\E[X_\ell(t)]\to\mu^*$ as $t\to\infty$ for some desired set-point $\mu^*>0$.
  \end{enumerate}
\end{objectives*}
The first requirement is fairly standard. Indeed, ergodicity is the analogue of having a globally attracting fixed point for deterministic dynamics (i.e. global stability) and is required here so that the closed-loop network is well-behaved, in the sense that it reaches stationarity \agedit{starting from any initial distribution}. The second requirement is more specific to stochastic processes, as even if the means converge, the variance can \agedit{still grow unboundedly with time}, which would mean that the actual dynamics of the process (its sample-paths) is not well-behaved, \red{rendering the controller of little practical utility}. Finally, the third statement encapsulates our desired objective of perfect adaptation (or set-point tracking), i.e. that the population mean \agedit{of the regulated species $\X{\ell}$} approaches a fixed homeostatic value $\mu^*$.

\subsection*{The Controller Reactions}

We propose the following controller network (Fig.~\ref{fig:cloud}) inspired from the deterministic networks proposed in \cite{Oishi:10b}:
\begin{equation}\label{eq:p:controller}
  \begin{array}{ccc}
    \underbrace{\phib\rarrow{\mu}\Z{1}}_{\textnormal{reference}},&\quad&\underbrace{\phib\rarrow{\theta X_\ell}\Z{2}}_{\textnormal{measurement}},\\ \underbrace{\Z{1}+\Z{2}\rarrow{\eta}\phib}_{\textnormal{comparison}}, &\quad& \underbrace{\phib\rarrow{kZ_1}\X{1}}_{\textnormal{actuation}}.
  \end{array}
\end{equation}
where $\X{\ell}$ is the \agedit{species \emph{measured} by the controller}, which is identical to the \emph{regulated species} in the current setup. The species $\Z{1}$ and $\Z{2}$ are referred to as the \emph{controller species}. 
Note that the topology of the controller network belongs to a family of four control topologies (see the supplementary material) depending on \agedit{respective roles of the two controller species}.


Although inspired from \cite{Oishi:10b}, the above network has a different philosophy. Besides the fact that the current setting is stochastic, the main difference lies in the way the network interacts with the environment. While the goal of \cite{Oishi:10b} was the biomolecular implementation of linear input-output systems, the goal here is to control a reaction network. In this regard, the birth-reactions of $\Z{1}$ and $\Z{2}$ clearly differ from the way they are defined in \cite{Oishi:10b}. We now clarify the role and meaning of each of the \agedit{controller} reactions:
\begin{enumerate}
  \item The first reaction is the \emph{reference reaction} (or set-point) which (partially) sets the value of the reference $\mu^*=\mu/\theta$. This value is implemented as the birth-rate of species $\Z{1}$.
  \item The second reaction is the \emph{measurement reaction} and takes the form of a pure-birth reaction with a rate proportional to the current population of the regulated species $\X{\ell}$\footnote{Note that it can also be implemented in terms of the catalytic reaction
  $\X{\ell}\rarrow{\theta}\X{\ell}+\Z{2}$.}. It is referred to as the \emph{measurement reaction} as the rate of increase of the population of $\Z{2}$ reflects the population \agedit{size} of $\X{\ell}$.
  \item The third reaction implements the \emph{comparison reaction}, \agedit{which decrements the molecular counts of $\Z{1}$ and $\Z{2}$ by one each. The rate constant for this reaction is} $\eta$ that can be tuned. The main role of this reaction is to correlate both the populations of $\Z{1}$ and $\Z{2}$ and to prevent them from growing without bounds. This reaction can be viewed as a \emph{compare and substract operation} since when both $\Z{1}$ and $\Z{2}$ have positive populations (comparison), then \agedit{this reaction decreases their respective population-sizes by one (subtraction)}, thereby preserving the difference \agedit{of their population-sizes}.
   \item The last reaction is the \emph{actuation reaction}, which implements the way the controller acts on the system, i.e. by acting on the birth-rate of the actuated species $\X{1}$\footnote{This can also be represented by the catalytic reaction  $\Z{1}\rarrow{k}\Z{1}+\X{1}$.}. The parameter $k$ is also a tuning parameter of the controller.
\end{enumerate}

The above controller has been chosen with an implementability constraint in mind, as it is expressed as plausible reactions that may be implemented in-vivo to perform \emph{in-vivo control}. It will be shown later that the proposed controller exhibits strong robustness properties which make its implementation much easier than other types of controllers that require the fine tuning of their reaction rates (see the supplementary material). In-vitro control is also possible using, for instance, DNA strand displacement  \cite{Chen:13}. In-silico control \cite{Khammash:11,Uhlendorf:12} can be considered as well, whenever the population size of regulated species $\X{\ell}$ can be measured in real-time from outside of the cell(s) using, for instance, microscopy.


\subsection*{Guaranteed Performance Properties of the Controlled Network}

\red{We now consider the dynamics of the \emph{closed-loop network}, which is formed by interconnecting the open-loop network $(\X{},\lambda,\zeta)$ with the controller \eqref{eq:p:controller}. We can represent the dynamics by a Markov process $\{  (X(t) , Z(t) ) : t \geq 0  \}$, where for each $t$, $X(t) =(X_1(t) ,\dots, X_d(t) )$ denotes the molecular counts of the $d$ network species and $Z(t) = (Z_1(t),Z_2(t) )$ denotes the molecular counts of the two controller species. Our control objective is to steer the first-order moment $\E[X_\ell(t)]$ corresponding to species $\X{\ell}$ to a desired set-point $\mu^*$. The dynamics of the first-order moments $\E[X(t)] = (\E[X_1(t)], \dots, \E[X_d(t)])$ and $\E[Z(t)] = (\E[Z_1(t)] , \E[Z_2(t) ])$ is described by the following system of ordinary differential equations (ODEs)
\begin{equation}\label{eq:1st_moment}
  \begin{array}{rcl}
    \frac{d \E[X(t)]}{d t}  &=& \sum_{k=1}^K  \zeta_k\E[ \lambda_k( X(t)   ) ] + k \E[ Z_1(t)] e_1, \\
      \frac{d \E[Z_1(t)]}{d t} & =&  \mu - \eta \E[Z_1(t)Z_2(t)],\\
\frac{d \E[Z_2(t)]}{d t} & = & \theta \E[X_\ell (t )] - \eta \E[Z_1(t)Z_2(t)],
  \end{array}
\end{equation}
%
%
%
where $e_1$ is the $d$-dimensional vector whose first component is $1$ and the rest are zero. Note that this system of ODEs is not \emph{closed} because there is no equation for the dynamics of $\E[Z_1(t)Z_2(t)]$. Moreover if the propensity functions $\lambda_k$'s are nonlinear functions of the state-variables, additional quantities whose dynamics is not captured by these equations will be encountered. Attempting to ``close" this system by adding equations for the dynamics of all these additional quantities and $\E[Z_1(t)Z_2(t)]$ will again lead to another set of quantities with unrepresented dynamics. This is known as the \emph{moment-closure problem} in the literature, a well-known barrier for the direct analysis and simulation of moment equations.}

\red{As stated, it is unclear why the proposed controller structure involves an integral action. To emphasize  this feature, let us define $\delta_Z(t):=Z_1(t)-Z_2(t)$. By subtracting the two last equations in \eqref{eq:1st_moment}, we get that
\begin{equation}\label{eq:deltaZ}
  \dfrac{d\E[\delta_Z(t)]}{dt}=\mu-\theta\E[X_\ell(t)].
\end{equation}
Defining the set-point \emph{tracking error} as $e(t):=\mu/\theta-\E[X_\ell(t)]$ and integrating the above expression over the interval $[0,t]$ yields
\begin{equation}\label{eq:integralaction}
  \E[\delta_Z(t)]=\theta\int_0^t e(s)ds+\E[\delta_Z(0)],
\end{equation}
which allows us to conclude that the first-order moment of the control input, $\E[u(t)]=k\E[Z_1(t)]$, will depend on the integral of the set-point tracking error $e(t)$ through the term $\E[\delta_Z(t)]$. Even though the control input only partially depends on the integrator state \eqref{eq:integralaction}, the stability of the closed-loop system proved in the main results indicates that the integral action is preserved.}

\red{Our approach in this paper is to find conditions ensuring that the Markov process $\{  (X(t) , Z(t) ) : t \geq 0  \}$ describing the dynamics of closed-loop reaction network is ergodic. This means that starting from any initial state, the probability distribution of the state $(X(t),Z(t))$ converges to a unique stationary distribution $\pi$ as $t \to \infty$. Under fairly general conditions, ergodicity also implies that first and second-order moments of the state $(X(t),Z(t))$ converge to their steady-state values, which can be computed by evaluating the expectation $\E_\pi$ with respect to the stationary distribution $\pi$.  Furthermore, as $t \to \infty$, the right hand side of \eqref{eq:deltaZ} tends to 0 and this yields the expression:
\begin{equation}\label{eq:equiI}
   \mu-\theta\E_\pi[X_\ell] = 0
\end{equation}
where $\E_{\pi} [X_\ell] = \lim_{ t \to \infty} \E[X_\ell(t)]$. From \eqref{eq:equiI}, we can immediately conclude that $\E_\pi[X_\ell] = \mu/\theta = \mu^*$, our desired objective.
This shows that if the closed-loop network dynamics is ergodic, then our proposed controller automatically imposes the set-point tracking property, $\E[X_\ell(t)]\to \mu^*$ as $t\to\infty$, regardless of the initial conditions. Obtaining such a behavior is the main rationale behind integral control. The controller we propose performs integral action on the dynamics, and achieves set-point tracking and perfect adaptation properties. Note that we demonstrate these properties of our controller without solving the first-order moment equations \eqref{eq:1st_moment}, thereby circumventing the moment-closure problem mentioned above.}

\red{The following result, proved in the supplementary material, establishes conditions under which a general stochastic biochemical reaction network can be controlled using the proposed controller network.}
\begin{theorem}[\bf General Network Case]\label{th:general_CL}
\red{Consider an open-loop reaction network $(\X{},\lambda,\zeta)$ and assume that for some given values of its parameters, and the parameters of the controller network, the closed-loop network, formed by augmenting the open-loop network with the controller reactions \eqref{eq:p:controller}, is ergodic and has uniformly bounded first- and second-order moments.} \red{Then, asymptotic set-point tracking is achieved, i.e. $\E[X_\ell(t)]\to\mu/\theta\ \textnormal{as}\ t\to\infty$.}
\end{theorem}

\red{The main challenge in the above result lies in verifying the ergodicity of a given closed-loop network dynamics.  In what follows, we provide simple conditions that allow us to check this property using efficient computational techniques, such as linear programming. This is done for two main classes of networks, namely those consisting of \emph{unimolecular} reactions and those consisting of \emph{bimolecular} reactions. In the unimolecular case, for each reaction $k$ the propensity function $\lambda_k(x)$ is an affine function of the state variable $x=(x_1,\dots,x_d)$. Hence we can express each $\lambda_k(x)$ as
\begin{align*}
\lambda_k(x) = \sum_{i=1}^d w_{ki} x_i+ w_{k0},
\end{align*}
where $w_{k0}$ is a nonnegative constant and $w_{ki}$'s are some real numbers for $i=1,\dots,d$. Define a $K \times d$ matrix $W$ with entries $w_{ki}$ and let $S$ be the $d\times K$ matrix whose $k$-th column is the stoichiometry vector $\zeta_k$ for reaction $k$. Also let $w_0$ be the $d$-dimensional vector whose $k$-th component is $w_{k0}$. Regarding each vector as a column-vector, for any state-vector $x$ we can write $$\sum_{k=1}^K \lambda_k(x) \zeta_k  =  S W x + Sw_0,$$
which allows us to express the first equation in \eqref{eq:1st_moment}, with $Z_1(t) \equiv 0$, as
    \begin{equation}\label{eq:linpos:nominalnew}
      \begin{array}{rcl}
        \dfrac{d\E[X(t)]}{d t}&=&SW\E[X(t)] +Sw_0.
      \end{array}
    \end{equation}
 This linear system of ODEs describes how the first-order moments of the open-loop network will evolve without the control action. If the matrix $SW$ only has eigenvalues with negative real parts (we say in this case that the matrix $SW$ is Hurwitz stable), then this system is \emph{asymptotically stable} and the first-order moment vector $\E[X(t)]$ converges as $t \to \infty$. To fulfil our control objective, it is necessary that this system be asymptotically stable. This is because our controller can only act positively on the open-loop network and hence it cannot stabilize an unstable system.}

\red{Recall that $\mathcal{S} \subset \N_0^d$ is the state-space for the Markovian dynamics of the open-loop reaction network. This state-space is irreducible if any state in $\mathcal{S}$ can be reached from any other state in $\mathcal{S}$ by a sequence of reactions having positive propensities at all the intermediate states. A simple example is the single-species birth-death process for which from any state value we can reach any larger state value by a sequence of birth reactions and any smaller state value by a sequence of death reactions. The irreducibility of the state-space is a necessary condition for ergodicity \cite{Meyn:93}, which holds for many reaction networks in the literature; see e.g. \cite{Briat:13i}. For a unimolecular open-loop network with an irreducible state-space $\mathcal{S}$, the asymptotic stability of the linear system \eqref{eq:linpos:nominalnew} is equivalent to the ergodicity of the open-loop network.}

\red{However in order to ensure the ergodicity of the closed-loop network, we need a couple of other conditions that are both very natural for our control problem. The first condition is that the open-loop system be \emph{output controllable} \cite{Ogata:70}, which simply means, in our case, that the molecular count of the regulated species $\X{\ell}$ responds to changes in the molecular count of the actuated species $\X{1}$. Mathematically, we can express this condition as
\begin{align}
\label{nonzeroimpluseresp}
[(SW)^{-1}]_{\ell 1}\ne0
\end{align}
where $[(SW)^{-1}]_{\ell 1}$ is the component at column $1$ and row $\ell$ of the inverse of matrix $SW$.
%
The second condition that we need is that the set-point $\mu^* = \mu/\theta$ should be accessible by the dynamics of $\E[X_\ell(t)]$. Since our controller can only act positively on the system \eqref{eq:linpos:nominalnew}, this accessibility condition can fail if some components of the input vector $Sw_0$ are \emph{too large} (see also the discussion in the supplementary material). Technically we can check this accessibility condition by ensuring that there exists a positive constant $c$ and a $d$-dimensional vector $v = (v_1,\dots,v_d)$ with positive entries such that each component of the vector $ v^T(SW+cI_d)$, where $I_d$ is the $d \times d$ identity matrix, is strictly negative and
   \begin{align}\label{eq:offset}
 \mu^* =  \dfrac{\mu}{\theta}>\dfrac{v^TSw_0}{c v_\ell}.
  \end{align}
As discussed in \cite{Briat:13i}, the above condition can be checked using linear programming techniques.}

\red{The following result, proved in the supplementary material, provides conditions under which a {\em unimolecular} open-loop reaction network can be controlled using the controller network \eqref{eq:p:controller}:}
\begin{theorem}[\bf Unimolecular Network Case] \label{th:p:affine:nominal}
\red{Suppose that the open-loop reaction network $(\X{},\lambda,\eta)$ is unimolecular and its state-space $\mathcal{S}$ is irreducible. Furthermore assume that the linear system \eqref{eq:linpos:nominalnew} of ODEs is asymptotically stable ($SW$ is Hurwitz stable) and output controllable (condition \eqref{nonzeroimpluseresp} holds), and that the desired set-point $\mu/\theta$ is accessible (condition \eqref{eq:offset} holds).}

\red{Then, for any $k,\eta>0$,  the closed-loop reaction network is ergodic and asymptotic set-point tracking is achieved, i.e. $\E[X_\ell(t)]\to\mu/\theta\ \textnormal{as}\ t\to\infty$. Moreover, the steady-state values of the first-order moments $\E[X(t)] = (\E[X_1(t)], \dots, \E[X_d(t)])$ and $\E[Z_1(t)]$ are
$$\lim_{t \to \infty}\E[X(t)]=\dfrac{\mu(SW)^{-1}e_1}{\theta [(SW)^{-1}]_{\ell 1}}\quad \textnormal{and}\quad \lim_{t \to \infty} \E[Z_1(t)]= \dfrac{-\mu}{\theta [(SW)^{-1}]_{\ell 1}}.$$}
\end{theorem}

A more general version of this result and its extension to a class of bimolecular networks are provided in the supplementary material. In light of Theorem \ref{th:p:affine:nominal}, several favorable properties for the controller and the closed-loop network can now be \red{highlighted and expounded}.

\textbf{Ergodicity, set-point tracking and bounded first- and second-order moments.} These are the main properties sought in \red{our statement of control objectives}. \red{Moreover, as stated in Theorem~\ref{th:p:affine:nominal}, the average population-sizes of species $\X{1},\dots,\X{d}$ and $\Z{1}$, at steady-state, are uniquely defined by the set-point $\mu/\theta$ and the parameters of the open-loop network, implying that these quantities are also regulated by our antithetic integral controller.}

\textbf{Robustness.} Robustness is a fundamental \red{requirement} which ensures that some properties for the closed-loop network are preserved, even in the presence of model uncertainties. This concept is critical in biology as the environment is fluctuating or \agedit{noisy} and \red{only poorly known models are typically} available. The obtained results can automatically guarantee the preservation of all the properties stated in Theorem~\ref{th:p:affine:nominal}, even in such constraining conditions.

\textbf{\red{Well-behaved single-cell tracking dynamics}.} Ergodicity ensures that the population average at stationarity is equal to the asymptotic value of the time-average of any single-cell trajectory; see e.g. \cite{Briat:13i}. We can therefore conclude that the proposed controller achieves two goals simultaneously, as it can ensure robust set-point tracking at both the population and single-cell levels. As a consequence, the controller will also ensure single-cell set-point tracking in the presence of cell events such as cell-division when certain conditions are met (see the supplementary material for more details).

\red{\textbf{Innocuousness of the controller.} An inaccurate implementation of controller parameters may sometimes lead to an unstable behavior for the closed-loop system. However, the fact that the conditions of Theorem \ref{th:p:affine:nominal} are independent of $k$ and $\eta$ tells us that the proposed controller will both preserve the ergodicity of the (possibly poorly known) open-loop network and ensure set-point tracking/adaptation regardless of the values of its parameters, provided that the open-loop network satisfies the conditions of Theorem \ref{th:p:affine:nominal}. This property is crucial in biology, as identifying models and implementing specific reaction rates (even approximately) are difficult tasks. This peculiar and non-standard property is referred here as \emph{innocuousness} and it is illustrated in Figure \ref{fig:det}, where the deterministic and stochastic dynamics of the same controlled networks are compared (see also the supplementary material).}

\red{\textbf{Low metabolic load.} \agedit{Even if the controller works in the low copy-number regime, it does not necessarily imply a low metabolic load for the cell. Indeed if there is fast creation and annihilation of the controller species $\Z{1}$ and $\Z{2}$, then it will result in many energy consuming futile cycles, which can impose a heavy metabolic burden on the cell, even though the dynamics is still in the low copy-number regime.} However, it can be shown (see the supplementary material) that the power consumption at stationarity of the controller reactions, denoted by $\bar{P}$, can be expressed in the case of unimolecular networks as
\begin{equation}
  \bar{P}=\mu\left(\alpha_1+\alpha_2+\alpha_3\right)+\dfrac{\mu}{\theta}\dfrac{\alpha_4}{|[(SW)^{-1}]_{\ell 1}|}
\end{equation}
where $\alpha_1,\alpha_2,\alpha_3$ and $\alpha_4$ are (positive) weights associated with the reference, measurement, comparison and actuation reactions, respectively, that represent the energy cost of each reaction. Interestingly, only the first three terms depend on $\mu$ while the last term depends on the ratio $\mu/\theta$ and some network parameters. This last term, however, would also be present in the case of the production of $\X{1}$ at a constitutive rate equal to $\mu/(\theta | e_\ell^T(SW)^{-1}e_1|)$ (which would lead to the same steady state value for the controlled output, but no adaptation properties). Hence, the effective metabolic load of our controller is equal to the first three terms and is only proportional to $\mu$.  A low metabolic load can therefore be easily achieved by first setting $\mu$ to a small value and then adjusting the set-point value with $\theta$.}

\textbf{Circumventing moment closure difficulties.} Finally, we emphasize that using the proposed approach, the moment closure problem does not arise, as the \agedit{main conclusions} (e.g. ergodicity, set-point tracking and robustness) directly follow from stochastic analysis tools and the structure of the controller, thereby avoiding altogether the framework of the moment equations (see the supplementary material).


\subsection*{Application to Gene Expression Control: set-point Tracking and Perfect Adaptation}

The goal of this example is to demonstrate that set-point tracking and perfect adaptation can be ensured with respect to any change in the parameters of the gene expression network
\begin{equation}\label{eq:gene_exp}
    \X{1}\rarrow{\gamma_1}\phib,\quad \X{1}\rarrow{k_2}\X{1}+\X{2},\quad \X{2}\rarrow{\gamma_2}\phib
\end{equation}
where $\X{1}$ denotes the mRNA and $\X{2}$ the measured/regulated species (see Fig.~\ref{fig:adaptation}). The following result is proved in the supplementary material:
\begin{proposition}
  For any positive values of the parameters $k,k_2,\gamma_1,\gamma_2,\eta,\theta$ and $\mu$, the controlled gene expression network \eqref{eq:p:controller}-\eqref{eq:gene_exp} is ergodic, has bounded and globally converging first- and second-order moments and
  \begin{equation}
    \E[X_2(t)]\to\dfrac{\mu}{\theta}\quad \textnormal{as}\quad t\to\infty.
  \end{equation}
\end{proposition}

\begin{figure}[H]
    \centering
  \includegraphics[width=0.9\textwidth]{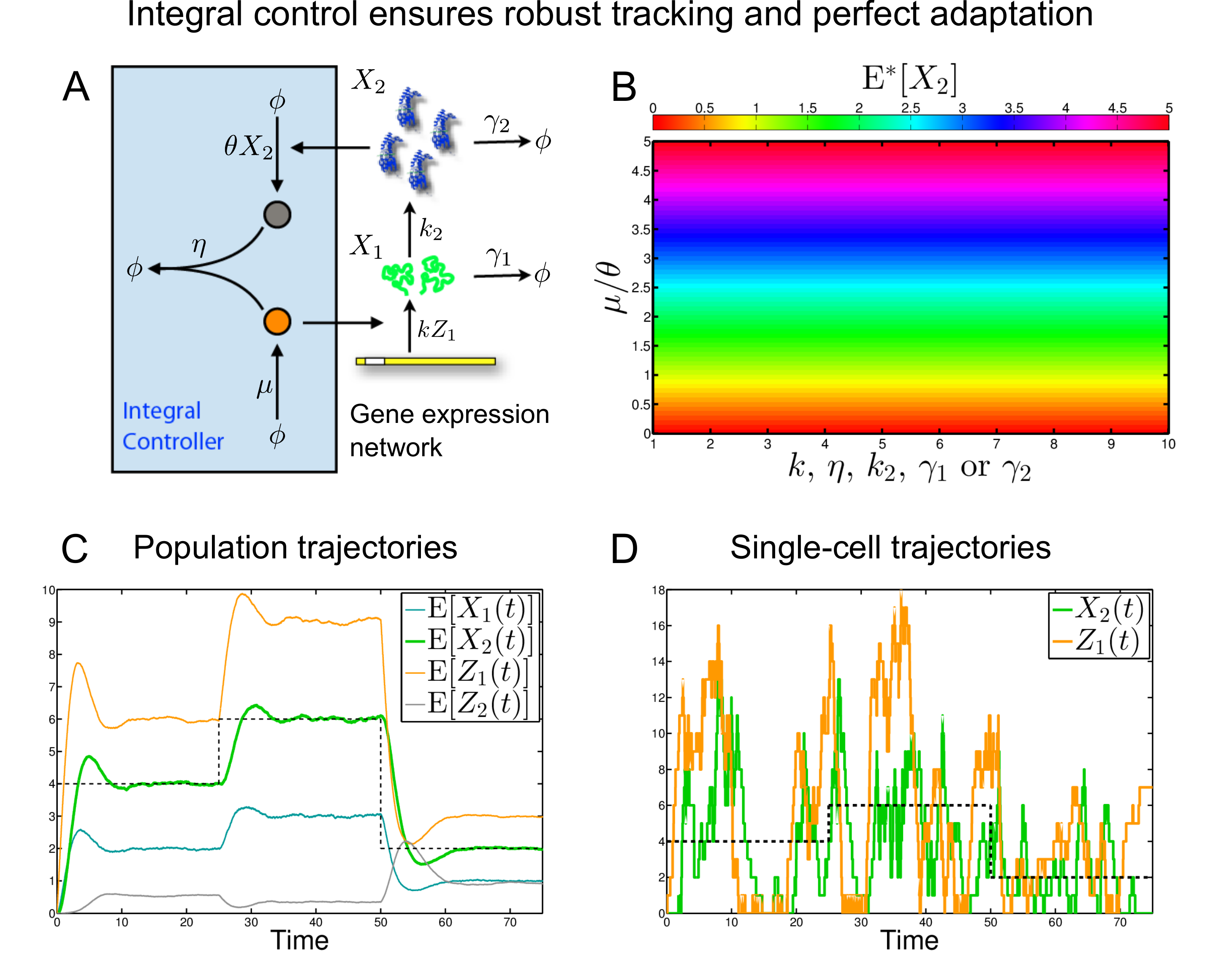}
  \caption{ \textbf{A.} The controlled gene expression network \eqref{eq:gene_exp} with the proposed antithetic integral controller \eqref{eq:p:controller}. \textbf{B.} The closed-loop reaction network shows perfect adaptation (at stationarity) with respect to any changes in the parameters of the network as we have that $\E_\pi[X_2]=\mu/\theta$ for any values of the parameters $k,\eta,k_2,\gamma_1$ and $\gamma_2$ where $\E_\pi[X_2]$ denotes the mean number of molecules of $\X{2}$ at stationarity. \textbf{C.} The controlled-output $\E[X_2(t)]$ of the closed-loop network tracks the reference value (in black-dash). The mean population of input species $\E[Z_1(t)]$ adapts automatically to changes in the reference value $\mu^*=\mu/\theta$ without requiring re-implementation.  \textbf{D.} Single-cell trajectories, although strongly affected by noise, still have an underlying regularity ensuring the convergence of the moments at the population level. All simulations have been performed using Gillespie's stochastic simulation algorithm with the parameters $k=1$, $\gamma_1=3$, $k_2=2$, $\gamma_2=1$, $\theta=1$ and $\eta=50$.}\label{fig:adaptation}
\end{figure}

\begin{figure}[H]
    \centering
   \includegraphics[width=0.9\textwidth]{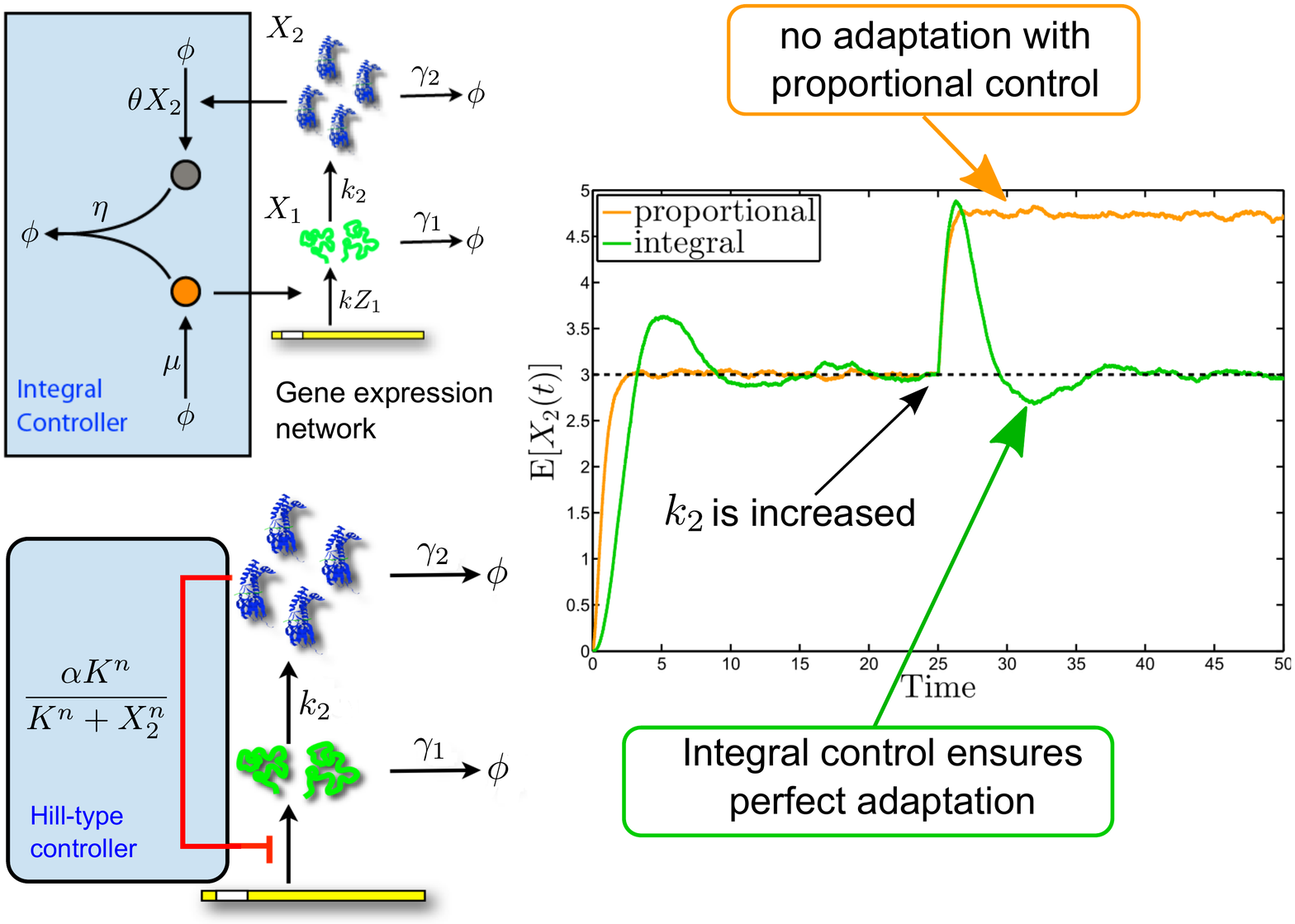}
  \caption{Comparison between the proposed antithetic integral controller \eqref{eq:p:controller} and the Hill-type static controller \eqref{eq:proportional_hill}. The simulation is performed using parameters initialized to $\mu=3$, $\theta=1$, $k=1$, $\gamma_1=3$, $k_2=3$, $\gamma_2=1$ and $\eta=50$ for the antithetic integral controller \eqref{eq:p:controller} and $n=1$, $\alpha=8.22$ and $K=3$ for the Hill-type static controller. The averaging is performed over 8000 cells simulated with Gillespie's stochastic simulation algorithm. At $t=25$s, the value of $k_2$ jumps from 3 to 6. While the proposed antithetic integral controller shows perfect adaptation, the Hill-type static controller is unable to return to the mean value of the population of $\X{2}$ before the stimulus. This demonstrates the advantage of the integral feedback strategy over the Hill-type strategy.}\label{fig:proptxt}
\end{figure}

\red{However, when we consider a Hill-type static control scheme of the form
\begin{equation}\label{eq:proportional_hill}
  \phib\rarrow{f(X_2)}\X{1}\quad \textnormal{with}\quad   f(X_2)=\frac{\alpha K^n}{K^n+X_2^n}
  \end{equation}
where $K,\alpha$ are positive parameters and $n$ is a positive integer, we obtain the results depicted in Fig.~\ref{fig:proptxt}. We can see that, as opposed to the antithetic stochastic integral controller, perfect adaptation is not ensured by the Hill-type static controller. Note that even though the comparison is only made for the gene expression network and this specific choice for the Hill-type static controller, it is a matter of fact that, in general, such controllers can not ensure perfect adaptation; see the supplementary material for some theoretical arguments and different control schemes.}

\subsection*{Noise as a Stabilizing Agent: Stochastic vs. Deterministic Population Control}

\begin{figure*}
    \centering
    \includegraphics[width=\textwidth]{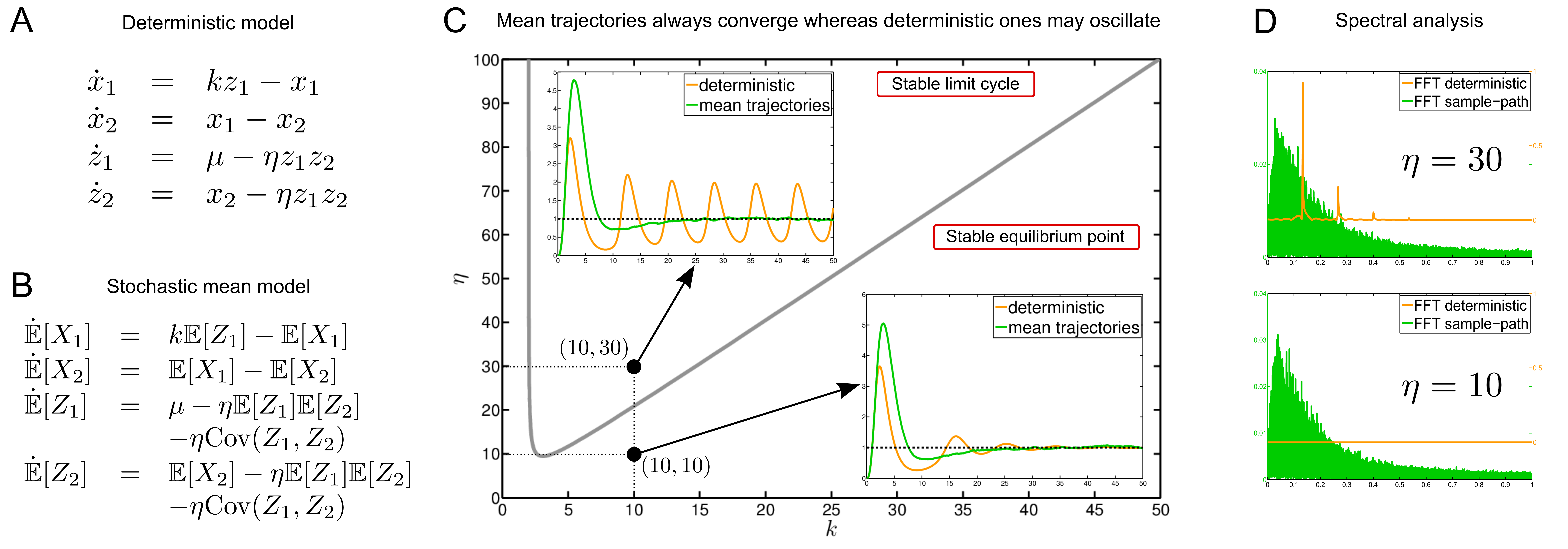}
 \caption{\textbf{A.} Deterministic model for the gene expression network \eqref{eq:gene_exp} with $k_2=\theta=\gamma_1=\gamma_2=1$. \textbf{B.}  Mean model for the gene expression network \eqref{eq:gene_exp} with $k_2=\theta=\gamma_1=\gamma_2=1$. \textbf{C.} The deterministic dynamics bifurcates from a unique stable equilibrium point when the controller parameters $(k,\eta)$ are chosen below the bifurcation curve into a stable limit-cycle when the controller parameters are chosen above. The first-order moments, however, always converge to the desired steady state value for the regulated species, here $\mu=1$, regardless of the values of the controller parameters. This can be explained by the presence of the stabilizing covariance term in the model for the stochastic means. \textbf{D.} While the frequency content at stationarity of the deterministic dynamics dramatically changes when crossing the bifurcation curve, the frequency content of the sample-paths remains qualitatively the same. In this regard, the controller can be considered to perform the same way in both cases. This demonstrates the superiority and the central role of the noise in the stabilizing properties of the proposed stochastic integral controller.}\label{fig:det}
\end{figure*}

\red{Here we demonstrate a striking effect of noise as an agent for dynamic stabilization at the population level. We do this by comparing here the results that we obtain} to those we would have obtained in the deterministic setting (see Fig.~\ref{fig:det}). To this aim, we again consider the gene expression network \eqref{eq:gene_exp} and we set $k_2=\gamma_1=\gamma_2=1$ for simplicity. We then get the deterministic model and stochastic mean model depicted in Fig.~\ref{fig:det}-\textbf{A} and Fig.~\ref{fig:det}-\textbf{B}, respectively. \red{It is important to emphasize that the deterministic model represents here the evolution of the mean concentration of the species over a population of deterministically behaving cells with identical initial concentrations. On the other hand, the stochastic mean model represents the evolution of the mean number of the species over a population of stochastically behaving cells with identical initial molecular counts. Note that, by virtue of the ergodicity property, the stochastic mean model will always converge to its unique steady-state value regardless of the different initial conditions for the individual cells. This property does not hold in general for a deterministic model representing the average concentrations of a population of deterministically oscillating cells.} The stochastic mean model has been obtained using the identity $\E[Z_1Z_2]=\E[Z_1]\E[Z_2]+\textnormal{Cov}(Z_1,Z_2)$ where the covariance term is nonzero as the random variables are not independent. If such a term would be zero, then we would recover the deterministic dynamics, but, due to noise, we can see in Fig.~\ref{fig:det}-\textbf{C} that while the deterministic dynamics may exhibit oscillations, the dynamics of the first-order moment is always globally converging to the desired steady-state value. As a final comment, we note that if we were closing the \agedit{moment equations} in Fig.~\ref{fig:det}-\textbf{C} by neglecting the second-order cumulant, then we would fail in predicting the correct behavior of the first-order moments.  This demonstrates the central role of the noise in the stabilizing properties of the proposed stochastic \agedit{antithetic} integral controller. \red{The noise can hence be viewed as here a stabilizing agent through the randomness it adds \agedit{to the dynamics}, allowing then for their systematic compensation at the population level, regardless of the initial conditions and the system parameters. This phenomenon is entirely \agedit{due to stochasticity} and it does not generically occur \agedit{in the deterministic setting}.}

\section*{Discussion.}

\red{A general control theory for stochastic biochemical reaction networks with tailored mathematical concepts and tools has been missing. We believe that a well-grounded \emph{biomolecular control theory}  would enable a deeper understanding of biological regulation at the molecular level and could pave the way for an efficient and systematic rational design of synthetic genetic circuits that function to regulate and steer cellular dynamics at the molecular level. For such circuits, we
propose the term ``cybergenetic", which combines the genetic nature of the system with its cybernetic function of dynamic steering and control.
In this article, we take a first step in the development of such a theory by addressing one of the central dynamic control motifs: integral feedback. The methods we developed are the product of a synthesis of ideas from control theory, probability theory, linear algebra and optimization theory. Even though our findings are specific to the class of integral controllers we consider, they may serve as the foundation on which more general \emph{biomolecular control theory} \agedit{can be developed} -- one that deals with a larger class of stochastic dynamic controllers and networks. Indeed, numerical experiments performed on more general networks lying outside the scope of the developed theory tend to support this claim (see the supplementary material).}

\begin{figure*}
    \centering
    \includegraphics[width=0.8\textwidth]{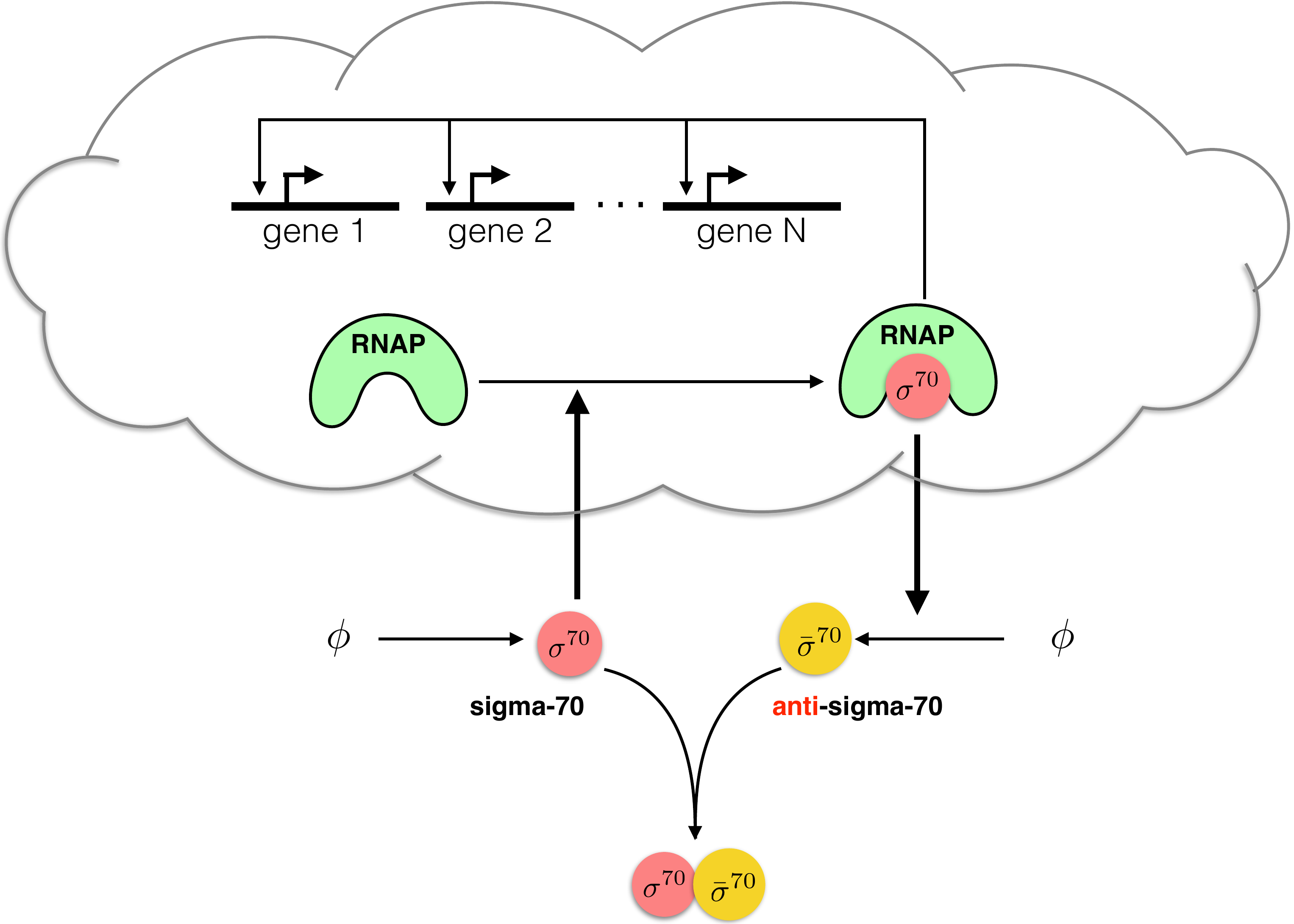}
 \caption{{\bf Regulation of the housekeeping genes in {\em E. coli}}. An endogenous circuit that may employ the control strategy proposed in this article is the $\boldsymbol{\sigma^{70}}$ regulation circuit. The sigma-factor $\boldsymbol{\sigma^{70}}$ binds RNA polymerase core in a complex that controls the expression of housekeeping genes during exponential growth conditions. The anti-$\boldsymbol{\sigma^{70}}$ factor, Rsd, whose gene is controlled by $\boldsymbol{\sigma^{70}}$, binds with a very strong affinity to $\boldsymbol{\sigma^{70}}$, sequestering it away from RNA polymerase. According to the theory put forth in this article, both interactions result in the tight regulation of the concentration of the $\boldsymbol{E\sigma^{70}}$ complex through negative feedback. While the average abundance of RNA polymerase is about 4600 molecules per cell \cite{Bakshi:12}, the average abundance of the complex $\boldsymbol{E\sigma^{70}}$ (regulated species) is approximately 700 molecules per cell \cite{Jishage:99}. At the same time, the abundance of \emph{free}-$\boldsymbol{\sigma^{70}}$ is computed to be fewer than one molecule per cell.} \label{sigma-70}
\end{figure*}

\red{Although the proposed control motif served as a frame around which the theoretical ideas of molecular stochastic control \agedit{were developed}, a cybergenetic circuit implementing such a motif may be of biological significance in its own right, a possibility which that we explore next for both endogenous and synthetic regulation.
Indeed, the simplicity of the mechanism and the remarkable regulation properties that it confers raises the question whether such a motif could have evolved for the purpose of endogenous regulation. The species $\Z{1}$ and $\Z{2}$ may be RNA or protein, but they must act to effectively ``annihilate" each other. Of course this annihilation need not be physical, as long as the two species act to render each other functionless. This could occur, for example, as a result of irreversible binding of $\Z{1}$ and $\Z{2}$, where the new complex effectively sequesters both from performing their function. Remarkably, one finds just such a possibility in the sigma-factor-mediated {\em E. coli} regulatory system. One example is the sigma factor $\boldsymbol{\sigma^{70}}$, which binds to RNA polymerase core enzyme ($\mathbf{E}$) enabling it to recognize the promoter of a host of housekeeping genes (see Figure~\ref{sigma-70}). It is well known that $\boldsymbol{\sigma^{70}}$ has a corresponding anti-sigma factor (Rsd) \cite{Trevino:13}, denoted here by $\boldsymbol{\bar{\sigma}^{70}}$, which binds $\boldsymbol{\sigma^{70}}$ very tightly, sequestering it away from {\bf E}.  Transcription of the Rsd gene is itself controlled by a $\boldsymbol{\sigma^{70}}$-dependent promoter \cite{Jishage:99}. This architecture fits very well our integral feedback regulation motif and strongly suggests that the core enzyme complex with $\boldsymbol{\sigma^{70}}$
($\boldsymbol{E \sigma^{70}}$) is the object of tight regulation. This in turn points to a corresponding regulation of the expression of all the housekeeping genes whose promoters are recognized by $\boldsymbol{\sigma^{70}}$. When binding reaction rates in the literature \cite{Maeda:00,Bakshi:12} are used for $\boldsymbol{\sigma^{70}}$ and RNAP along with estimated average total numbers, one comes up with a very low ($<1$) average copy number of {\em free}-$\boldsymbol{\sigma^{70}}$ per cell, which would indicate that the system operates in the stochastic regime. Intriguingly, the anti-$\boldsymbol{\sigma^{70}}$ factor is known to be involved in the transition to stationary phase, during which time its concentrations are significantly elevated, leading to the down-regulation of the housekeeping genes. One way to achieve this is through the increase of the rate of transcription of anti-sigma factor Rsd during transition to stationary phase, leading to the decrease in the set-point for $\boldsymbol{E\sigma^{70}}$ (in our analysis this corresponds to the quantity $\mu/\theta$).  Supporting this hypothesis is the fact that transcription of the $\boldsymbol{\sigma^{70}}$ gene is known to be controlled by a $\boldsymbol{\sigma^{38}}$-dependent promoter \cite{Jishage:99}, where $\boldsymbol{\sigma^{38}}$ is recognized as the master regulator for adaptation to stationary phase transcription.  While the full regulatory details of $\boldsymbol{\sigma^{70}}$ remain unknown, it is quite likely that $\boldsymbol{\sigma^{70}}$ and $\boldsymbol{\bar{\sigma}^{70}}$ interaction as described here plays a central part in that regulation.  It is also quite possible that such a regulatory motif is not uncommon in biology, as several anti-sigma factors have been found in a number of bacteria, including {\em E. coli} and Salmonella, as well as in the T4 bacteriophage.}

Beyond endogenous circuits, our cybergenetic motif presents opportunities for applications in synthetic biology.  Until now most of the synthetic regulatory circuits have relied on proportional action--a control scheme that fails to ensure perfect adaptation in many practical situations. Moreover, existing theoretical studies of synthetic biological circuits mainly considered the deterministic setting, and hence they implicitly assumed large molecular abundances. However, the implementation of control circuits that rely on high component abundances severely impinges on the host circuit's material and energy resources, leading to increased metabolic burden which can affect both function and viability. Fortunately, this is largely avoidable, as effective control involves mostly information processing, which in principle requires little energy and material resource consumption. \red{The novel regulatory motif that we propose exhibits characteristics that provably ensure robust stability, robust set-point tracking, and robust perfect adaptation for the controlled network and is achieved with a low metabolic cost and with molecular species that can have very low abundances.} It can be used for both single-cell set-point tracking (on average) and for population control. Thanks to the innocuousness of the controller, it does not need to be fine tuned, and can therefore be used in many practical situations, e.g. when the controlled network is very poorly known. In this regard, the proposed controller maintains clear implementability advantages over controllers requiring parameter tuning. This latter property emerges from the random nature of the reactions, as its deterministic counterpart leads to oscillating trajectories when the controller parameters are located in a certain instability region. \red{In spite of this, this controller can still be used in a deterministic setting even though some of the properties, such the innocuousness property, are lost.} With this in mind, the proposed controller may find several applications within synthetic biology. An immediate one is the optimization of drug or fuel production in bioreactors; see e.g. \cite{Dunlop:10}. Currently simple control strategies, such as proportional feedback or constitutive production, are used in these applications. By utilizing slightly more complex controllers, such as the one proposed here,  dramatic improvements in the production process can be expected, thanks to their enhanced robustness properties.  Another important application example is the design of insulators; see e.g. \cite{Mishra:14}. It has indeed been shown that loading effects are often detrimental to modular design. Insulators are therefore needed in order to preserve function modularity. The proposed controller can be used as a buffering element in order to drive the output of a module to the input of another one. It can also be used as a constant signal generator that can be used to act on a network to be analyzed. The amplitude can be tuned by acting on the reference, which can be modified from outside the cell using light-induced techniques \cite{Khammash:11}.


The proposed controller, however, may have some drawbacks, as it seems to introduce some additional variance to the controlled process. Even though this extra variance is not detrimental to the current control objectives, it may be a problem if \agedit{the} goal is to reduce the variance over a cell population. Whether the variance can be reduced via a more optimal choice of parameters or through \agedit{some additional} controller reactions, remains a question for further research. In the end, some \agedit{``extra"} variance \agedit{due to the controller} may be unavoidable, as fundamental limitations to variance reduction with feedback \cite{Lestas:10} are likely to hold  for any molecular control circuit. \agedit{Even so, controller noise should not detract from the tremendous promise of designing novel stochastic control circuits at the molecular level, where the dynamic properties, benefits, and limitations seem to be exquisitely different from those at the macroscopic scale.}

\section*{Acknowledgments}
The authors are grateful to Adam Arkin for pointing out the similarity of the tight sigma/anti-sigma factor binding reactions to our feedback controller reactions. The authors acknowledge funding support from the Human Frontier Science Program Grant RGP0061/2011 and the Swiss National Science Foundation grant 200021-157129.


\section*{Conflict of Interest}

No conflicts of interest.

\bibliographystyle{plain}

\end{document}